\documentclass[11pt]{amsart}
\numberwithin{equation}{section}
\newtheorem{theorem}{Theorem}
\newtheorem{definition}{Definition}

\numberwithin{theorem}{section} \numberwithin{lemma}{section}
 \numberwithin{definition}{section}
\numberwithin{proposition}{section}

\def\al{\aligned}
\def\eal{\endaligned}
\def\M{{\bf M}}
\def\be{\begin{equation}}
\def\ee{\end{equation}}
\def\lab{\label}

\def\M{\bf M}

\def\al{\aligned}
\numberwithin{equation}{section}

\begin{document}

\tracingpages 1
\title[volume ratio]{\bf bounds on volume growth of geodesic balls under Ricci flow}
\author{Qi S. Zhang}
\address{Department of
Mathematics, University of California, Riverside, CA 92521, USA }
\date{July 2011}

\begin{abstract}
We prove a so called $\kappa$ non-inflating property for Ricci
flow, which provides an upper bound for volume ratio of geodesic
balls over Euclidean ones, under an upper bound for scalar
curvature. This result can be regarded as the opposite statement
of Perelman's $\kappa$ non-collapsing property for Ricci flow.
These two results together imply volume doubling property for
Ricci flow without assuming Ricci curvature lower bound.

\end{abstract}
\maketitle
\section{Statement of result and proof}

In \cite{P:1}, Perelman proved the fundamental $\kappa$
non-collapsing property for Ricci flow. One version of it roughly
says that the volume ratio between a geodesic ball and Euclidean
ball with the same radius is bounded from below by a positive
constant, provided that the scalar curvature is bounded from above
in a space time cube.

In this short note, we prove that the opposite result is also
true. i.e. the above volume ratio is bounded from above by a
positive constant, provided that the scalar curvature is bounded
from above in a space time cube. In the case of normalized Ricci
flow on compact K\"ahler manifolds with positive first Chern
class, the upper bound holds for all time. An upper bound for the
volume ratio is useful in the study of K\"ahler Ricci flow. See
for example the papers \cite{Se:1}, \cite{CW:1} and the references
therein.  The current result seems to remove one of the obstacles
in the program to prove convergence results, although many other
obstacles remain.

To make the statement precise, let's introduce notations and definition.
We use $\M$ to denote a compact Riemann manifold and $g(t)$ to denote the
metric at time $t$; $d(x, y, t)$ is the geodesic distance under $g(t)$;
$B(x, r, t) = \{ y \in {\M} \ | \ d(x, y, t) < r \}$  is the geodesic ball of radius $r$, under metric $g(t)$,
centered at $x$, and $|B(x, r, t)|_{g(t)}$ is the volume of
$B(x, r, t)$ under $g(t)$;  $d\mu_{g(t)}(x)$ is the volume element. We also
reserve $R=R(x, t)$ as the scalar curvature under $g(t)$.

\begin{definition}
 A smooth, compact, $n$ dimensional Ricci flow  $({\M}, g(t))$ is called $\kappa$ non-inflated at the point
 $(x_0, t_0)$ under scale $\rho$ if the following statement holds.

 For any $r \in (0, \rho)$, suppose:

 1. the Ricci flow is defined in the space time cube
 \[
 \{ (x, t) \, | \, d(x, x_0, t_0) < r, \quad t \in [t_0-r^2, t_0] \},
 \]

 2. for some positive constant $\alpha$,
 $R(x, t) \le \frac{\alpha}{t_0-t}$ for all $(x, t)$ in the above cube.

Then, there exists a positive constant $\kappa$, which may depend on $\alpha$ such that
\[
|B(x_0, r, t_0)|_{g(t_0)} \le \kappa r^n.
\]
 \end{definition}
\medskip

{\it \noindent Remark.  Recall that in the $\kappa$ non-collapsing property, the condition on the scalar curvature is
$R(x, t) \le \frac{1}{r^2}$ in the space time cube. Obviously this condition is included in our condition
$R(x, t) \le \frac{\alpha}{t_0-t}$ in the same space time cube.
}
\medskip

The main result of the note is the following theorem. Even though the proof is very short,
it actually uses a combination of several results: a variation of
global bounds for the fundamental solution of
the conjugate heat equation (\cite{CZ:1}) whose proof relies on Perelman's Harnack inequality in \cite{P:1} and uniform Sobolev inequality under Ricci flow, and also Perelman's scalar
curvature and diameter bound for K\"ahler Ricci flow, and the general idea that heat kernel lower bound implies volume upper bound in \cite{GHL:1}.

\begin{theorem}
\lab{thm1.1}
(a).  Let $({\M}, g(t))$,  $\partial_t g_{ij} = - 2 R_{ij}$, $t \in [0, t_0]$ be a smooth, compact, $n$ dimensional Ricci flow.
Then for any $x_0 \in \M$, the Ricci flow is $\kappa$ non-inflated at $(x_0, t_0)$ under
scale $\sqrt{t_0}$. Here $\kappa$ depends only on  $g(0)$, $t_0$ and the constant $\alpha$
in the bound for the scalar curvature.

(b). Let $({\M}, g(t))$,  $\partial_t g_{ij} = -  R_{ij} + g_{ij}$, $t \in [0, \infty)$ be a smooth, compact, $n$ real dimensional, normalized Ricci flow on compact K\"ahler manifolds with positive first
Chern class. There exists a positive constant $\kappa>0$, which depends only on
the initial metric $g(0)$ such that
\[
|B(x, r, t)| \le \kappa r^n
\]for all $x \in \M$, $r>0$ and $t>0$.
\end{theorem}

\proof (of part (a)).

\noindent{Step 1.}
\medskip

Picking any $r \in (0, \sqrt{t_0})$, we assume:

1.  the Ricci flow is defined in the space time cube
 \[
 Q(x_0, t_0, r)= \{ (x, t) \, | \, d(x, x_0, t_0) < r, \quad t \in [t_0-r^2, t_0] \},
 \]

 2. for some positive constant $\alpha$,
 $R(x, t) \le \frac{\alpha}{t_0-t}$ for all $(x, t)$ in the above cube $Q(x_0, t_0, r)$.

Let $l$, $t$ be two moments in time such that $0<l<t \le t_0$, and $x, z \in \M$.
Let $G=G(z, l; x, t)$ be the fundamental solution of the conjugate heat equation
\be
\lab{conjheat}
\Delta u -R u + \partial_l u=0
\ee which is coupled with the Ricci flow. Fixing $z, l$, we know that
$G$, as a function of $x, t$ satisfies the heat equation i.e., for $t>l$,
\[
\Delta_x G(z, l; x, t) - \partial_t G(z, l; x, t) = 0.
\]Hence
\[
\al
&\frac{d}{dt} \int_{\M} G(z, l; x, t) d\mu_{g(t)}(x) =
\int_{\M} [\Delta_x G(z, l; x, t) -R(x, t) G(z, l; x, t)] d\mu_{g(t)}(x) \\
&= -  \int_{\M} R(x, t) G(z, l; x, t) d\mu_{g(t)}(x).
\eal
\]From the scalar curvature equation
\[
\Delta R + 2 |Ric|^2 - \partial_t R=0,
\]we deduce
\[
\Delta R + \frac{2}{n} R^2 - \partial_t R \le 0,
\]which implies, via the maximum principle
that, either $R(\cdot, 0) \ge 0$ or
\[
\min R(\cdot, t) \ge \frac{1}{\left( 1/\min R(\cdot, 0) \right) - (2 t/n)}.
\]Therefore, either $\frac{d}{dt} \int_{\M} G(z, l; x, t) d\mu_{g(t)}(x) \le 0$ or
\[
\frac{d}{dt} \int_{\M} G(z, l; x, t) d\mu_{g(t)}(x)
 \le \frac{1}{\left( -1/\min R(\cdot, 0) \right) + (2 t/n)} \int_{\M} G(z, l; x, t) d\mu_{g(t)}(x),
 \]which yields
 \be
 \lab{intG}
 \int_{\M} G(z, l; x, t) d\mu_{g(t)}(x) \le 1+C (1+t-l)^{n/2}.
 \ee Here $C$ only depends on $\min R(\cdot, 0)$ and $n$, and $C=0$ when $R \ge 0$.

 \medskip
 \noindent{Step 2.}

 Next we prove the following heat kernel bounds on $G(z, l; x, t)$ which is similar to
 Theorem 2.1 in \cite{CZ:1}.  The method is also similar. The improvement is on the coefficients of the bounds which rely on
  the initial metric $g(0)$ instead of on $g(l)$.  This will be useful in proving the volume ratio bound.
 \be
 \lab{Gxiajie}
  \frac{c_1 J(t)}{ (t-l)^{n/2}}
e^{-2 c_2 \frac{d(z, x, t)^2}{t-l}}  e^{- \frac{1}{ \sqrt{t-l}}
\int^{t}_{l} \sqrt{t-s}
 R(x, s) ds} \le G(z, l; x, t) \le \frac{c^{-1}_1 J^{-1}(t)}{ (t-l)^{n/2}}.
 \ee Here
 \[
 J=J(s)=\exp[-\alpha - s \beta- s \, \sup R^-(\cdot, 0)],
 \]and $\alpha$ and
  $\beta$ are positive constants depending only on the Sobolev constants of $({\M}, g(0))$ and the
  infimum of Perelman's $F$ entropy for $({\M}, g(0))$. The proof of this theorem uses uniform Sobolev inequality under
 Ricci flow and Perelman's differential Harnack inequality for $G$ (\cite{P:1}). Moreover, if
 the scalar curvature is positive, then $J(t)$ is independent of $t$.

Recall that as a function of $(x, t)$, $G=G(z, l; x, t)$ is the
fundamental solution of the forward heat equation associated with
the Ricci flow, i.e.,

\begin{equation}\label{forheat}
\begin{cases}\begin{array}{lll}\frac{\partial}{\partial
t}g(t)&=&-2Ric,\\
u_t&=&\Delta u.\end{array}\end{cases}\end{equation}
We study the forward heat equation (\ref{forheat}) first.\\

Let $u=u(x,t)$ be a positive solution to (\ref{forheat}). Given $T
>0$ and $t \in (l, T)$, define
\[
p(t)=(T-l)/(T-t),
\] so $p(l)=1$ and $p(T)=\infty$. By direct
computation, using the idea of Davies,
\[
\al
\partial_t \Vert u \Vert_{p(t)} &= \partial_t \bigg[ \bigg(\int_{\bf M}
u^{p(t)}(x, t) d\mu_{g(t)} \bigg)^{1/p(t)} \bigg]\\
&=-\frac{p'(t)}{p^2(t)} \Vert u \Vert_{p(t)} \ln \int_{\bf M}
u^{p(t)}(x, t) d\mu_{g(t)}
    +\frac{1}{p(t)} \bigg( \int_{\bf M} u^{p(t)}(x, t) d\mu_{g(t)} \bigg)^{(1/p(t))-1}
          \\
          &\qquad \times \bigg[ \int_{\bf M} u^{p(t)} (\ln u) p'(t) d\mu_{g(t)} +
     \int_{\bf M} u^{p(t)-1} ( p(t) \Delta u - R u ) d\mu_{g(t)} \bigg].
\eal
\] Using integration by parts on the
term containing $\Delta u$ and multiplying both sides by $p^2(t)
\Vert u \Vert^{p(t)}_{p(t)}$, we infer
\[
\al &p^2(t)  \Vert u \Vert^{p(t)}_{p(t)}
\partial_t \Vert u \Vert_{p(t)}\\
=& -p'(t) \Vert u \Vert^{p(t)+1}_{p(t)} \ln \int_{\bf M}
u^{p(t)}(x, t) d\mu_{g(t)}
  + p(t) \Vert u \Vert_{p(t)} p'(t) \int_{\bf M}
      u^{p(t)} \ln u (x, t) d\mu_{g(t)}\\
      &-p^2(t)(p(t)-1) \Vert u \Vert_{p(t)} \int_{\bf M}
            u^{p(t)-2} |\nabla u|^2(x, t) d\mu_{g(t)}
              - p(t) \Vert u \Vert_{p(t)} \int_{\bf M}  R(x, t)
            u^{p(t)}(x, t) d\mu_{g(t)}.
\eal
\]Dividing both sides by $ \Vert u \Vert_{p(t)}$, we arrive at
\[
\al & p^2(t)  \Vert u \Vert^{p(t)}_{p(t)}
\partial_t \ln \Vert u \Vert_{p(t)}\\
=& -p'(t) \Vert u \Vert^{p(t)}_{p(t)} \ln \int_{\bf M} u^{p(t)}
d\mu_{g(t)}
  + p(t)  p'(t) \int_{\bf M}
      u^{p(t)} \ln u  d\mu_{g(t)}\\
      &-4[p(t)-1]  \int_{\bf M}
            |\nabla (u^{p(t)/2}) |^2 d\mu_{g(t)}
              - p(t)  \int_{\bf M}  R
            (u^{p(t)/2})^2 d\mu_{g(t)}.
\eal
\]
Define $v(x,t) = \frac{u^{p(t)/2}}{\Vert u^{p(t)/2} \Vert_2}$,
then $\Vert v \Vert_2 = 1$ and $$v^2\ln v^2=p(t) v^2\ln u-2v^2 \ln
\Vert u^{p(t)/2} \Vert_2.
$$ Merging the first two terms on the
righthand side of the above equality and dividing both sides by
$\Vert u \Vert^{p(t)}_{p(t)}$, we find that
\[
\al & p^2(t)
\partial_t \ln \Vert u \Vert_{p(t)}\\
=& p'(t)   \int_{\bf M} v^2 \ln v^2  d\mu_{g(t)}
         -4(p(t)-1)  \int_{\bf M}
            |\nabla v |^2 d\mu_{g(t)}
              - p(t)  \int_{\bf M}  R v^2 d\mu_{g(t)}\\
      =& p'(t)   \int_{\bf M} v^2 \ln v^2 d\mu_{g(t)}
         -4[p(t)-1]  \int_{\bf M}
            (|\nabla v |^2 +\frac{1}{4} R v^2) d\mu_{g(t)} - \int_{\bf M}  R v^2 d\mu_{g(t)}.
 \eal
\]Notice the following relations,
\[
\frac{4(p(t)-1)}{p'(t)} = \frac{4 (t-l) (T-l-(t-l))}{T-l} \le T-l, \qquad  \frac{1}{p'(t)}=\frac{(T-l)^2}{T}
\le T,
\]
Hence
\[
\al & p^2(t)
\partial_t \ln \Vert u \Vert_{p(t)} \\
\le& p'(t) \bigg[  \int_{\bf M} v^2 \ln v^2  d\mu_{g(t)} -
\frac{4(p(t)-1)}{p'(t)} \int_{\bf M}
            (|\nabla v |^2 +\frac{1}{4} R v^2) d\mu_{g(t)}
+  T \sup R^-(x, t) \bigg]. \eal
\]Take $\epsilon$ such that
\[
\epsilon^2= \frac{4(p(t)-1)}{p'(t)} \le T-l \] in the
log-Sobolev inequality (6.2.8) in Section 6.2 of \cite{Zshu:1}, we deduce that
\[
p^2(t)
\partial_t \ln \Vert u \Vert_{p(t)} \le
 p'(t) \bigg[ -n \ln \sqrt{4 (p(t)-1)/p'(t)} + L(t) +  T \sup
 R^-(x, 0) \bigg],
 \]where, due to $\epsilon^2 \le T-l \le T$,
 \[
 \al
 L(t)  &\doteq (t +\epsilon^2) \beta  + \alpha\\
&\le 2 T  \beta + \alpha \\
& \doteq L(T), \eal
 \] for some positive constants $\alpha=\alpha(A_0, B_0, \lambda_0, n)$ and
 $\beta=\beta(A_0, B_0, \lambda_0, n)$. Here $A_0$, $B_0$ are the coefficients in the standard Sobolev
 inequality for $({\M}, g(0))$ and $\lambda_0$ is the infimum of Perelman's $F$ entropy for $({\M}, g(0))$.
 We stress that these depend only on initial metric.
 Here we also used the fact that
 \[
 \sup
 R^-(x, t) \le \sup
 R^-(x, 0),
 \]which is a consequence of maximum principle and the
 evolution equation of
 scalar curvature
  $\partial_t R=\Delta R + 2 |Ric|^2$.

Observe that $p'(t)/p^2(t)=1/(T-l)$ and
\[
4(p(t)-1)/p'(t)=4(t-l) [T-l-(t-l)]/(T-l).
\]
 Hence we have
 \[
\partial_t \ln \Vert u \Vert_{p(t)} \le
 \frac{1}{T-l} \bigg\{ - \frac{n}{2} \ln [4 (t-l) [T-l-(t-l)]/(T-l)] + L(T) + T \sup
 R^-(x, 0) \bigg\}.
 \]This implies, after integrating from $t=l$ to $t=T$, that
 \[
 \ln  \frac{\Vert u(\cdot, T) \Vert_\infty}{\Vert u(\cdot, l)
 \Vert_1}
 \le - \frac{n}{2} \ln (4 (T-l)) + L(T) + T \sup
 R^-(x, 0)+n.
 \]Since
 \[
 u(x, T) = \int_{\bf M}  G(z, l; x, T) u(z, l) d\mu_{g(l)},
 \]
the above inequality implies that
 \begin{equation}
 \lab{ondiagbound}
 G(z, l, x, T) \le \frac{\exp[L(T)+  T \sup
 R^-(x, 0)]}{(4  (T-l))^{n/2}},
 \end{equation}
 where $L(T)$ is defined above as
 \[
 L(T) = 2  T \beta+\alpha.
\] Since $T$ is arbitrary, we get the upper bound by defining $J^{-1}=L(T)+  T \sup
 R^-(x, 0)$.  Note the constants  $\beta$ and $\alpha$ may have changed by a factor.

In case $R(x, 0)>0$, by Section 6.2 in \cite{Zshu:1},
we have $\beta=0$.
 So the above bound becomes
\begin{equation}
\label{ondiagbound2}
 G(z, l; x, T) \le \frac{\exp(\alpha )}{(4 \pi (T-l))^{n/2}},
 \end{equation}proving the upper bound.\\

Next we prove the lower bound. Let $t<t_0$  and $u=u(x, t) \equiv
G(x, t; x_0, t_0)$.   We claim that for
 a constant $C>0$,
 \[
 G(x_0, t; x_0, t_0) \ge  \frac{C}{\tau^{n/2}} e^{- \frac{1}{2 \sqrt{\tau}} \int^{t_0}_{t} \sqrt{t_0-s}
 R(x_0, s) ds} .
 \]where $\tau = t_0-t$ here and later in the proof.
  To prove this inequality,  define  a function $f$ by
 \[
 (4 \pi \tau)^{-n/2} e^{-f} = u.
 \]We need to apply Perelman's differential Harnack inequality for the
 fundamental
 solution along any smooth space-time curve $\gamma (t)$ (see \cite[Corollary 9.4]{P:1}).
 Here we pick the curve
 $\gamma(t)$ to be the fixed point $x_0$, we have,
 \[
 -\partial_t f(x_0, t) \le \frac{1}{2} R(x_0, t) - \frac{1}{2 \tau} f(x_0, t).
 \]For any $t_2<t_1<t_0$, we integrate the above inequality to get
 \[
 f(x_0, t_2) \sqrt{t_0-t_2} \le  f(x_0, t_1) \sqrt{t_0-t_1}  + \frac{1}{2} \int^{t_1}_{t_2} \sqrt{t_0-s}
 R(x_0, s) ds.
 \]When $t_1$ approaches $t_0$,  $f(x_0, t_1)$ stays bounded since
$G(x_0, t_1; x_0, t_0) (t_0-t_1)^{n/2} $ is bounded between two
positive constants, which is a direct consequence of the standard asymptotic formula for $G$
(for example, see \cite[Chapter 24]{chowetc3}).  Hence for any $t \le t_0$, we
have
 \[
 f(x_0, t)  \le  \frac{1}{2 \sqrt{t_0-t}} \int^{t_0}_{t} \sqrt{t_0-s}
 R(x_0, s) ds.
 \]Consequently
\be
\lab{Gx0x0}
G(x_0, t; x_0, t_0) \ge \frac{c}{(4 \pi \tau)^{n/2}} e^{-
\frac{1}{2 \sqrt{t_0-t}} \int^{t_0}_{t} \sqrt{t_0-s}
 R(x_0, s) ds} .
 \ee

As observed earlier, the function $G(x_0, t; \cdot, \cdot)$ is a
solution to the standard heat equation coupled with Ricci flow,
which is the conjugate of the conjugate heat equation. i.e.,
 \[
 \Delta_z G(x,
t; z; l) - \partial_l G(x, t; z, l) =0,
\] here $\Delta_z$ is with respect to the metric $g(l)$. Therefore
it follows from \cite[Theorem 3.3]{z06} or \cite[Theorem
5.1]{ch2009} that, for $\delta>0, c_1, c_2>0$, and $y_0 \in \M$,
\[
G(x_0, t; x_0, t_0) \le c_1 G^{1/(1+\delta)} (x_0, t, y_0, t_0)
K^{\delta/(1+\delta)}  e^{c_2 d^2(x_0, y_0, t_0)/\tau},
\]where $K =\sup_{M \times [(t_0+t)/2, t_0]} G(x_0, t, \cdot, \cdot)$.
The upper bound
\[
K \le \frac{c J^{-1}(t_0)}{(t_0-t)^{n/2}},
\]
 together with the lower  bound (\ref{Gx0x0}) imply that, with $\delta = 1$,
\[
G(x_0, t; y_0, t_0) \ge c_1 \frac{J(t_0)}{ (t_0-t)^{n/2}}
e^{-2 c_2 d(x_0, y_0, t_0)^2/\tau}  e^{- \frac{1}{ \sqrt{t_0-t}}
\int^{t_0}_{t} \sqrt{t_0-s}
 R(x_0, s) ds},
\] which is the desired  lower bound.

\medskip
\noindent{Step 3.}
\medskip

 In (\ref{Gxiajie}), we take $z=x_0$, $t=t_0$ and $l=t_0-r^2$ where $r$ is
 the given number in $(0, \sqrt{t_0})$. By the assumption on the scalar curvature $R(x, t)$ in the
 definition of $\kappa$ non-inflating,,
 we obtain, for $x$ such that $d(x_0, x, t_0) \le r$,
 \[
 \al
 G(x_0, t_0-r^2; x, t_0) & \ge  \frac{c_1 J(t_0)}{r^n}
e^{-2 c_2 }  e^{- \frac{1}{r}
\int^{t_0}_{t_0 - r^2} \sqrt{t_0-s}
 R(x, s) ds}\\
 &\ge \frac{c_1 J(t_0)}{r^n}
e^{-2 c_2 }  e^{- \frac{1}{r}
\int^{t_0}_{t_0 - r^2} \sqrt{t_0-s}
 \frac{\alpha}{t_0-s} ds}.
\eal
\] Thus, when $d(x_0, x, t_0) \le r$, we have
\[
G(x_0, t_0-r^2; x, t_0)  \ge \frac{c_1 J(t_0)}{r^n} e^{-2 c_2 - 2
\alpha}.
\]Substituting this to (\ref{intG}), we deduce
\[
\al
1+C (1+r^2)^{n/2}& \ge \int_{\M} G(x_0, t_0-r^2; x, t_0) d\mu_{g(t_0)}(x)\\
&\ge \int_{d(x_0, x, t_0) \le r} G(x_0, t_0-r^2; x, t_0) d\mu_{g(t_0)}(x) \\
& \ge \frac{c_1 J(t_0)}{r^n} e^{-2 c_2 - 2 \alpha} \int_{d(x_0, x,
t_0) \le r}  d\mu_{g(t_0)}(x).
 \eal
\]This implies
\[
| B(x_0, r, t_0)|_{g(t_0)} r^{-n} \le [1+ C (1+t_0)^{n/2}] e^{2 c_2 + 2 \alpha}  c^{-1}_1 J^{-1}(t_0).
\]Taking
\[
\kappa = [1+C (1+t_0)^{n/2} ]e^{2 c_2 + 2 \alpha}  c^{-1}_1 J^{-1}(t_0),
\] we
obtain
\[
| B(x_0, r, t_0)|_{g(t_0)} \le \kappa \, r^n
\]
proving part (a) of the theorem. Note $\kappa$ depends only on $t_0$ and $g(0)$ in general, and
if $R \ge 0$, then $\kappa$ only depends on $g(0)$, due to the aforementioned property on the constant $J$ in
(\ref{Gxiajie}), and the fact that in the expression of $\kappa$, the constant $C=0$ when $R \ge 0$.
\medskip

\proof (of part (b)).
Since the normalized Ricci flow is smooth, we only need to prove the result for $t \ge c_0$ for some
positive constant $c_0$. According to Perelman (see \cite{ST:1}), the scalar curvature $R=R(x, t)$ and diameter of
the manifold are uniformly bounded for all time. Using the scaling
\[
t=-\ln(1- 2 \tilde t), \quad g(t)=\frac{1}{1-2 \tilde t} {\tilde  g} (\tilde t),
\]we see that $ {\tilde  g} (\tilde t)$ is the standard Ricci flow in the time interval  $\tilde t \in [0, 1/2)$ such that
\[
{\tilde R}(x, \tilde t) \le \frac{\alpha}{1- 2 \tilde t},
\]where $\alpha$ is a positive constant. Pick $\tilde t \in [1/4, 1/2)$ and $\tilde r \in (0, 1/2)$. Then for all
$s \in [\tilde t-\tilde r^2, \tilde t]$ and $x \in \M$, we have
\[
{\tilde R}(x, s) \le \frac{\alpha/2}{(1/2)-s} \le \frac{\alpha}{\tilde t-s}.
\]
Now we can just apply part (a) of the theorem to conclude
\[
| \{ y \ | \ d(y, x, \tilde g(\tilde t)) < \tilde r \} |_{\tilde g(\tilde t)} \le \kappa \tilde r^n
\]where $\kappa$ depends only on the initial metric $g(0)$ and $\alpha$. This is so because the total length of time interval is $1/2$ for $\tilde t$. After scaling we obtain, for $r =
(1-2 \tilde t)^{-1/2} \tilde r=e^{t/2} \tilde r$,
\[
|B(x, r, t)|_{g(t)} \le \kappa r^n.
\]
Since $\tilde r$ can be any number in $(0, 1/2)$,
 we conclude that for $r \le e^{t/2}/2$, and all $t  \ge \ln 2$,
\[
|B(x, r, t)|_{g(t)} r^{-n} \le \kappa.
\]Since the diameter of $({\M}, g(t))$ is uniformly bounded, the above holds for all $r>0$ with perhaps a
different constant $\kappa$.
\qed

\medskip
{\bf Acknowledgement.}  We wish to thank Professors Xiaodong Cao
and Xiouxiong Chen for very helpful suggestions. After the paper
is circulated and posted on the arxiv, Professor Xiouxiong Chen kindly informed us that he
and Professor Bing Wang also obtained a similar result when the scalar curvature has an extra lower bound, which is just posted \cite{CW:2}. Thanks also go to Professor Hongliang Shao for checking the 
paper carefully and correcting a typo.

\bigskip

\noindent e-mail:  qizhang@math.ucr.edu

\enddocument